\documentclass[12pt,draft]{amsart}

\usepackage[cp1251]{inputenc}
\usepackage[T2A]{fontenc}
\usepackage[russian,english]{babel}
\usepackage{amsmath,amsfonts,amssymb}
\usepackage{geometry}
\newtheorem{theorem}{Theorem}

\newtheorem{lemma}{Lemma}

\theoremstyle{remark}

\theoremstyle{definition}

\title[General forms of the Menshov--Rademacher, Orlicz, and Tandori theorems]
{General forms of \\ the Menshov--Rademacher, Orlicz, and Tandori theorems \\ on
orthogonal series}

\author[V. A. Mikhailets, A. A. Murach]{Vladimir A. Mikhailets, Aleksandr A. Murach}

\address{Institute of Mathematics, National Academy of Sciences of Ukraine,
3 Tereshchenkivs'ka, Kyiv, 01601, Ukraine}

\email{mikhailets@imath.kiev.ua, murach@imath.kiev.ua}



\begin{document}

\maketitle

\begin{abstract} We prove that the classical Menshov--Rademacher, Orlicz, and
Tandori theorems remain true for orthogonal series given in the direct integrals of
measurable collections of Hilbert spaces. In particular, these theorems are true for
the spaces $L_{2}(X,d\mu;H)$ of vector-valued functions, where $(X,\mu)$ is an
arbitrary measure space, and $H$ is a real or complex Hilbert space of an arbitrary
dimension.
\end{abstract}

\section{Introduction}\label{sec1}

The Menshov--Rademacher theorem \cite{Menschoff23, Rademacher22} plays an important
role in the theory of orthogonal series. It states that the sequence
$(\log^{2}_{2}n)$ is a Weyl multiplier for convergence, almost everywhere
(a.e.) with respect to the Lebesgue measure, of a series in an arbitrary orthonormal system (ONS)
of real-valued functions given on a finite interval of the real axis. There are some
various theorems on unconditional convergence of orthogonal series. These results refine the
Menshov--Rademacher theorem (see, e.g., \cite[Ch.~2, \S~5]{Alexits61} and
\cite[Ch.~8, \S~2]{KashinSaakyan89}), where the Orlicz theorem \cite{Orlicz27}
occupies a special place. It gives a sufficient condition for the sequence
$(\omega_{n}\log^{2}_{2}n)$ to be a Weyl multiplier for the unconditional
convergence a.e. The Menshov--Rademacher and the Orlicz theorems are best possible
in the sense that their conditions cannot be weakened.

It is known (see, e.g., \cite{MoriczTandori96, Meaney07}) that the
Menshov--Rademacher theorem remains valid for series with respect to ONSs of
real-valued or complex-value functions given on an arbitrary measure space. This
also true \cite{11UMJ10} for the Orlicz theorem and for another known result on
unconditional convergence, the Tandori theorem \cite{Tandori62}.

The question arises whether these and others theorems on convergence of orthogonal
series are true in a more general setting of series with respect to ONSs of
vector-valued functions given on a measure space and taking values in a collection
of Hilbert spaces.

In the present paper, we will give a positive answer to this question for the
classical Menshov--Rademacher, Orlicz, and Tandori theorems.

Note that, in the case of orthogonal series in (complex-valued)
eigenfunctions of a self-adjoint elliptic operator defined on a closed compact
manifold $X$, the conditions of the Menshov--Rademacher and the Orlicz theorems and that
the function being expanded  belongs to the isotropic H\"ormander
spaces $H^{\psi}(X)$ are equivalent, where $\psi(t)=\log^{*}t$ or $\psi(t)=\varphi(t)\log^{*}t$,
respectively; see \cite{08MFAT1, 09Dop3} and \cite[Sec.~2.3.2]{MikhailetsMurach10}.
Here $\log^{*}t:=\max\{1,\log_{2}t\}$, whereas $\varphi(t)$, $t\geq1$, is a positive
increasing function that varies regularly at $+\infty$ in the sense of Karamata and
satisfies the condition
$$
\int\limits_{2}^{\infty}\,\frac{dt}{t\,(\log_{2}t)\,\varphi^{2}(t)}<\infty.
$$

\section{Statements of the main results}\label{sec2}

Let $X$ be an arbitrary measurable space with some $\sigma$-additive measure
$\mu\geq0$. The measure is not assumed to be finite or $\sigma$-finite. Let
$\{H(x):x\in X\}$ be a $\mu$-measurable collection of either real or complex Hilbert
spaces. This means that the function $\dim H(x)$, $x\in X$, takes only finitely or
countably many values (that are cardinal numbers) and that all the sets
$$
\{\,x\in X:\,\dim H(x)=\mathrm{const}\,\}
$$
are $\mu$-measurable. We consider the direct integral
$$
\mathbf{L}_{2}:=\int\limits_{X}^{\oplus}H(x)\,d\mu(x)
$$
of the $\mu$-measurable collection $\{H(x):x\in X\}$ (see, e.g., \cite[Ch.~7,
Sec.~1]{BirmanSolomjak87} and \cite[Ch.~2]{Nielsen80}). The space $\mathbf{L}_{2}$
is endowed with the inner product
$$
(f(\cdot),g(\cdot))_{2}:=\int\limits_{X}(f(x),g(x))_{H(x)}\,d\mu(x),
$$
which induces the norm $\|\cdot\|_{2}$.

If $H(x)\equiv H=\mathrm{const}$, then
$$
\mathbf{L}_{2}=L_{2}(X,d\mu;H)=L_{2}(X,d\mu)\otimes H.
$$
Thus, in this case, the space $\mathbf{L}_{2}$ consists of all classes of
$\mu$-equivalent vector-valued functions $f:X\rightarrow H$ that are strongly
measurable with respect to $\mu$ \cite[Ch.~V, Sec. 4]{Yosida65} and that
$$
\|f\|_{2}=\biggl(\,\int\limits_{X}\|f(x)\|_{H}^{2}\,d\mu(x)\biggr)^{1/2}<\infty.
$$

Let an ONS of vector-valued functions $\Phi:=(\varphi_{n})_{n=1}^{\infty}$ be
arbitrarily chosen in the space $\mathbf{L}_{2}$. We investigate the $\mu$-almost
everywhere ($\mu$-a.e.) convergence on $X$ of the orthogonal series
\begin{equation}\label{f1}
\sum_{n=1}^{\infty}\,a_{n}\,\varphi_{n}(x).
\end{equation}
Here all coefficients $a_{n}$ are either complex or real numbers; this depends on
whether all the spaces $H(x)$, $x\in X$ are complex or real. We set
$a:=(a_{n})_{n=1}^{\infty}$. Given $x\in X$, the convergence of the series
\eqref{f1} is regarded in the norm of $H(x)$.

Consider the majorant of partial sums of this series:
\begin{equation}\label{f2}
S^{*}(\Phi,a,x):=
\sup_{m\in\mathbb{N}}\,\bigl\|\,\sum_{n=1}^{m}\,a_{n}\,\varphi_{n}(x)\bigr\|_{H(x)},
\quad x\in X.
\end{equation}

Let us formulate the main results of the paper.

\begin{theorem}[a general form of the Menshov--Rademacher theorem]\label{th1}
Let a sequence of numbers $(a_{n})_{n=1}^{\infty}$ satisfy the condition
\begin{equation}\label{f3}
L:=\sum_{n=1}^{\infty}\,|a_{n}|^{2}\,\log_{2}^{2}(n+1)<\infty.
\end{equation}
Then the series \eqref{f1} converges $\mu$-a.e. on $X$, and moreover
\begin{equation}\label{f4}
\|S^{*}(\Phi,a,\cdot)\|_{2}\leq K\,\sqrt{L}.
\end{equation}
Here $K$ is a certain universal positive constant, one may take $K=4$.
\end{theorem}

This theorem was proved independently by D.~E.~Menshov \cite{Menschoff23} and
H.~Rademacher \cite{Rademacher22} in the case where
\begin{equation}\label{case}
X=(\alpha,\beta),\;\;-\infty<\alpha<\beta<\infty,\;\;\;\mu\,\;\mbox{is the Lebesgue
measure},\;\;\;H(x)\equiv\mathbb{R}.
\end{equation}
An exposition of their results are given, e.g., in G.~Alexits' \cite[Sec.
2.3.2]{Alexits61} and B.~S.~Kashin and A.~A.~Saakyan's \cite[Ch.~8, \S~1]{KashinSaakyan89} books. Note that the measures $\mu$ that are absolutely continuous
with respect to the Lebesgue measure are also allowed in \cite{Alexits61}. As it has
been mentioned, the Menshov--Rademacher theorem remains true for the ONSs of
real-valued or complex-valued functions given on an arbitrary measure space $X$.
Remark that a complete characterization of the sequences $(a_{n})_{n=1}^{\infty}$
such that the series \eqref{f1} converges a.e. for an arbitrary ONS in $L_{2}(X,
d\mu;\mathbb{R})$ is given by A.~Paszkiewicz \cite{Paszkiewicz10}.

The Menshov--Rademacher theorem is precise. In the situation \eqref{case},
D.~E.~Menshov \cite{Menschoff23} constructed an example of ONS
$(\varphi_{n})_{n=1}^{\infty}$ such that for every sequence of numbers
$(\omega_{n})_{n=1}^{\infty}$ satisfying
$$
1=\omega_{1}\leq\omega_{2}\leq\omega_{3}\leq\ldots,\quad
\lim_{n\rightarrow\infty}\,\frac{\omega_{n}}{\log_{2}^{2}n}=0
$$
there exists an a.e. divergent series of the form \eqref{f1} whose coefficients meet
the condition
$$
\sum_{n=1}^{\infty}\;|a_{n}|^{2}\,\omega_{n}<\infty.
$$
This result is presented, e.g., in the books \cite[Sec. 2.4.1]{Alexits61} and
\cite[Ch.~8, \S~1]{KashinSaakyan89} mentioned above.

Recall that the series \eqref{f1} is called \emph{unconditionally} convergent
$\mu$-a.e. on $X$ if the series
\begin{equation}\label{f5}
\sum_{n=1}^{\infty}\;a_{\sigma(n)}\,\varphi_{\sigma(n)}(x)
\end{equation}
converges $\mu$-a.e. on $X$ for an arbitrary permutation
$\sigma=(\sigma(n))_{n=1}^{\infty}$ of the set $\mathbb{N}$ of all positive
integers. Here the zero measure set of the points at which the series \eqref{f5}
diverges can depend  on the permutation~$\sigma$.

\begin{theorem}[a general form of the Tandori theorem]\label{th2}
Let a sequence of numbers $(a_{n})_{n=1}^{\infty}$ satisfy the condition
\begin{equation}\label{f6}
\sum_{k=0}^{\infty}\,\biggl(\;\sum_{n=\nu_{k}+1}^{\nu_{k+1}}\,
|a_{n}|^{2}\,\log_{2}^{2}n\biggr)^{1/2}<\infty,
\end{equation}
where $\nu_{k}:=2^{2^{k}}$. Then the series \eqref{f1} converges unconditionally
$\mu$-a.e. on $X$.
\end{theorem}

This theorem was proved by K.~Tandori \cite{Tandori62} in the situation
\eqref{case}. He also showed that his theorem is best possible in the following
sense. Given a (nonstrictly) decreasing sequence of positives numbers
$(a_{n})_{n=1}^{\infty}$, the series \eqref{f1} converges unconditionally a.e. for
each ONS $(\varphi_{n})_{n=1}^{\infty}$ in $L_{2}((0;1),dx,\mathbb{R})$ if and only
if \eqref{f6} holds. These K.~Tandori's results are presented in the book
\cite{KashinSaakyan89} (see Ch.~8, \S~2 and the remarks to Ch.~8).

A sufficient condition for the unconditional convergence of the series \eqref{f1}
can be expressed in the terms of the Weyl multipliers.

\begin{theorem}[a general form of the Orlicz theorem]\label{th3}
Let a sequence of numbers $(a_{n})_{n=1}^{\infty}$ and a (nonstrictly) increasing
sequence of positives numbers $(\omega_{n})_{n=1}^{\infty}$ satisfy the following
conditions:
\begin{gather}\label{f7}
\sum_{n=2}^{\infty}\;|a_{n}|^{2}\,(\log_{2}^{2}n)\,\omega_{n}<\infty,\\
\sum_{n=2}^{\infty}\;\frac{1}{n\,(\log_{2}n)\,\omega_{n}}<\infty. \label{f8}
\end{gather}
Then the series \eqref{f1} converges unconditionally $\mu$-a.e. on $X$.
\end{theorem}

Under the assumption \eqref{case}, Theorem~\ref{th3} is an equivalent formulation of
the Orlicz theorem \cite{Orlicz27}, which was suggested by P.~L.~Ulj'anov \cite[\S~4,
Sec.~1]{Uljanov63} (also see \cite[\S~9, Sec.~1]{Uljanov64}). The Orlicz theorem and
its proof can be founded, e.g., in G.~Alexits' book \cite[Sec. 2.5.1]{Alexits61}. As
K.~Tandori proved \cite{Tandori62}, this theorem is best possible in the sense that
the condition \eqref{f8} on the sequence $(\omega_{n})_{n=1}^{\infty}$ cannot be
weakened.

Note that both Theorems \ref{th2} and \ref{th3} remain true for each ONS of
complex-valued functions given on an arbitrary measure space $X$ \cite{11UMJ10}.

Theorems \ref{th1}, \ref{sec2}, and \ref{th3} will be proved in Sections \ref{sec4},
\ref{sec5}, and \ref{sec6}, resp. When proving Theorems \ref{th1} and \ref{th2}, we
use the classical scheme of argumentation set forth in \cite[Ch.~8, \S 1,
2]{KashinSaakyan89} for the case \eqref{case}. Theorem~\ref{th3} will be deduced
from Theorem~\ref{th2}. Previously, in Section~\ref{sec3} we establish a general
form of the Menshov--Rademacher inequality that plays a decisive role in the
proofs of Theorems \ref{th1} and \ref{sec2}.

\section{Menshov-Rademacher inequality}\label{sec3}

The proofs of Theorem~\ref{th1} and \ref{th2} are based on the following fact.

\begin{lemma}\label{lem1}

Let an integer $N\geq1$, finite ONS of vector-valued functions
$\Psi:=(\psi_{n})_{n=1}^{N}$ in $\mathbf{L}_{2}$, and a finite collection of numbers
$b:=(b_{n})_{n=1}^{N}$ be given arbitrarily. Then the function
\begin{equation}\label{eq9}
S^{*}_{N}(\Psi,b,x):=\max_{1\leq j\leq
N}\,\bigl\|\,\sum_{n=1}^{j}\,b_{n}\,\psi_{n}(x)\,\bigr\|_{H(x)},\quad x\in X,
\end{equation}
satisfies the inequality
\begin{equation}\label{eq10}
\|S^{*}_{N}(\Psi,b,\cdot)\|_{2}\leq(2+\log_{2}N)
\,\biggl(\;\sum_{n=1}^{N}\,|b_{n}|^{2}\biggr)^{1/2}.
\end{equation}
\end{lemma}

In the classical case \eqref{case}, the inequality \eqref{eq10} was obtained
independently by D.~E.~Menshov \cite{Menschoff23} and G.~Rademacher
\cite{Rademacher22} and then used by them in the proof of Theorem~\ref{th1} (see, e.g.,
the books \cite[Sec. 2.3.1, 2.3.2]{Alexits61} and \cite[Ch.~9,
\S~1]{KashinSaakyan89}). On the right-hand side of \eqref{eq10}, the factor
$C\log_{2}(N+1)$ with some universal constant~$C$ is used usually instead of
$2+\log_{2}N$. Note that this inequality is known for ONSs of real-valued or
complex-valued functions given on an arbitrary measure space $X$ (see, e.g.,
\cite[Theorem~3]{Moricz76} and \cite[Proposition 2.1]{Meaney07}).

\medskip

\emph{Proof of Lemma}~\ref{lem1}\emph{.} First we consider the case when $N=2^{r}$
for some integer $r\geq\nobreak1$. The general situation is easily reduced to this
case; this will be shown at the end of the proof.

Given an arbitrary number $j\in\{1,\,2,\ldots,2^{r}\}$, consider its binary
representation
\begin{equation*}
j=\sum_{k=0}^{r}\,\varepsilon_{k}\,2^{r-k},\quad\mbox{where}\;\;
\varepsilon_{k}:=\varepsilon_{k}(j)\in\{0,\,1\}.
\end{equation*}
Then every sum $\sum_{n=1}^{j}h_{n}$ of vectors in a real or complex Hilbert space
$H$ can be represented in the form
\begin{equation*}
\sum_{n=1}^{j}h_{n}=\sum_{k\,:\,\varepsilon_{k}\neq0}\;\;
\sum_{\sum\limits_{s=0}^{k-1}\varepsilon_{s}2^{r-s}
<n\leq\sum\limits_{s=0}^{k}\varepsilon_{s}2^{r-s}}\!\!h_{n}.
\end{equation*}
Whence, using the triangle inequality for the norm in $H$ and the Cauchy inequality
(both being applied to the external sum of $\leq r+1$ terms), we get:
\begin{equation*}
\begin{aligned}
\bigl\|\,\sum_{n=1}^{j}\,h_{n}\,\bigr\|_{H}&=
\bigl\|\,\sum_{k\,:\,\varepsilon_{k}\neq0}\;\;1\cdot
\sum_{\sum\limits_{s=0}^{k-1}\varepsilon_{s}2^{r-s}
<n\leq\sum\limits_{s=0}^{k}\varepsilon_{s}2^{r-s}}\!\!h_{n}\,\bigr\|_{H}\\
&\leq\sum_{k\,:\,\varepsilon_{k}\neq0}\;\;1\cdot
\bigl\|\,\sum_{\sum\limits_{s=0}^{k-1}\varepsilon_{s}2^{r-s}
<n\leq\sum\limits_{s=0}^{k}\varepsilon_{s}2^{r-s}}\!\!h_{n}\,\bigr\|_{H}\\
&\leq(r+1)^{1/2}\;\biggl(\:\sum_{k\,:\,\varepsilon_{k}\neq0}\;\bigl\|\,
\sum_{\sum\limits_{s=0}^{k-1}\varepsilon_{s}2^{r-s}
<n\leq\sum\limits_{s=0}^{k}\varepsilon_{s}2^{r-s}}
\!\!h_{n}\,\bigr\|_{H}^{2}\biggr)^{1/2}\\
&\leq(r+1)^{1/2}\;\biggl(\:\sum_{k=0}^{r}\,\sum_{p=0}^{2^{k}-1}\;
\bigl\|\,\sum_{n=p2^{r-k}+1}^{(p+1)2^{r-k}}h_{n}\,\bigr\|_{H}^{2}\biggr)^{1/2}.
\end{aligned}
\end{equation*}
Thus
\begin{equation}\label{f13}
\bigl\|\,\sum_{n=1}^{j}\,h_{n}\,\bigr\|_{H}^{2}\leq
(r+1)\;\sum_{k=0}^{r}\,\sum_{p=0}^{2^{k}-1}\;
\bigl\|\,\sum_{n=p2^{r-k}+1}^{(p+1)2^{r-k}}h_{n}\,\bigr\|_{H}^{2}.
\end{equation}

We apply this inequality to estimate the function \eqref{eq9}, which is represented
in the form
\begin{equation*}
S^{*}_{N}(\Psi,b,x)=
\bigl\|\,\sum_{n=1}^{j(x)}\,b_{n}\,\psi_{n}(x)\,\bigr\|_{H(x)},\quad x\in X;
\end{equation*}
here the number $j(x)\in\{1,\,2,\ldots,2^{r}\}$ is properly chosen for every fixed
$x\in X$. Setting $h_{n}:=b_{n}\psi_{n}(x)$ in \eqref{f13}, write:
\begin{equation*}
(S^{*}_{N}(\Psi,b,x))^{2}\leq(r+1)\;\sum_{k=0}^{r}\,\sum_{p=0}^{2^{k}-1}\;
\bigl\|\,\sum_{n=p2^{r-k}+1}^{(p+1)2^{r-k}}b_{n}\psi_{n}(x)\,\bigr\|_{H(x)}^{2},\quad
x\in X.
\end{equation*}

Integrating the latter inequality and using that $(\psi_{n})_{n=1}^{2^{r}}$ is an
ONS in $\mathbf{L}_{2}$, we have:
\begin{equation*}
\begin{aligned}
\|S^{*}_{N}(\Psi,b,\cdot)\|_{2}^{2}&\leq
(r+1)\;\sum_{k=0}^{r}\,\sum_{p=0}^{2^{k}-1}\;\int\limits_{X}\,
\bigl\|\,\sum_{n=p2^{r-k}+1}^{(p+1)2^{r-k}}b_{n}\psi_{n}(x)\,\bigr\|_{H(x)}^{2}\,d\mu(x)\\
&=(r+1)\;\sum_{k=0}^{r}\,\sum_{p=0}^{2^{k}-1}\;
\sum_{n=p2^{r-k}+1}^{(p+1)2^{r-k}}|b_{n}|^{2}=
(r+1)\;\sum_{k=0}^{r}\,\sum_{n=1}^{2^{r}}\,|b_{n}|^{2}=
(r+1)^{2}\,\sum_{n=1}^{2^{r}}\,|b_{n}|^{2}.
\end{aligned}
\end{equation*}
Thus
\begin{equation}\label{eq19}
\|S^{*}_{N}(\Psi,b,\cdot)\|_{2}^{2}\leq(r+1)^{2}\,\sum_{n=1}^{2^{r}}\,|b_{n}|^{2}.
\end{equation}
This, in view of $N=2^{r}$, yields the required estimate~\eqref{eq10}.

Now consider the general situation, when $N\geq1$ is an arbitrary integer. If $N=1$,
then Lemma~\ref{lem1} is trivial. Let $N\geq2$; then there exists an integer
$r\geq1$ such that $2^{r-1}<N\leq2^{r}$. Putting $a_{n}:=0$ for $N<n\leq2^{r}$, we
arrive at the above case, when the collection $(a_{n})$ consists of $2^{r}$ numbers.
Therefore, \eqref{eq19} holds with $r-1<\log_{2}N$; i.e., the required
inequality~\eqref{eq10} is fulfilled in the general situation.

Lemma~\ref{lem1} is proved.

\section{Proof of Theorem~$\ref{th1}$}\label{sec4}

Beforehand let us make a useful remark. Without loss of generality we may assume
that the measure $\mu$ is $\sigma$-finite. Indeed, since $\|\varphi_{n}\|_{2}=1$ for
each $n\geq1$, it follows that every set $\{x\in X:\|\varphi_{n}(x)\|_{H(x)}>1/j\}$,
with $j\in\mathbb{N}$, has a finite measure. Hence, $\mu$ is a $\sigma$-finite
measure on the set of all points $x\in X$ such that $\varphi_{n}(x)\neq0$ for at
least one index $n$. Outside this set all terms of the series \eqref{f1} are
zero-vectors. Therefore our assumption does not lead to any loss of generality in
the proofs.

Now let us show that the sequence
\begin{equation}\label{eq23}
S_{2^{k}}(x):=\sum_{n=1}^{2^{k}}\,a_{n}\,\varphi_{n}(x),\quad k=1,\,2,\,3,\ldots,
\end{equation}
converges for $\mu$-a.e. $x\in X$, and then we estimate the norm in
$L_{2}(X,d\mu;\mathbb{R})$ of the function
\begin{equation*}
S^{\star}(x):=\sup_{0\leq k<\infty}\,\|S_{2^{k}}(x)\|_{H(x)},\quad x\in X.
\end{equation*}

Let
\begin{equation*}
\chi_{k}(x):=\sum_{n=2^{k}}^{2^{k+1}-1}\,a_{n}\,\varphi_{n}(x),\quad x\in X,\quad
k=0,\,1,\,2,\,3,\ldots\,.
\end{equation*}
Since $(\varphi_{n})_{n=1}^{\infty}$ is an ONS in $\mathbf{L}_{2}$, we may write
\begin{equation*}
\|\chi_{k}\|_{2}^{2}=\sum_{n=2^{k}}^{2^{k+1}-1}\,|a_{n}|^{2}.
\end{equation*}
Hence, by the condition \eqref{f3}, we have
\begin{equation*}
\begin{aligned}
\sum_{k=0}^{\infty}\|\chi_{k}\|_{2}^{2}\,(k+1)^{2}&=
\sum_{k=0}^{\infty}(k+1)^{2}\sum_{n=2^{k}}^{2^{k+1}-1}\,|a_{n}|^{2}\\
&\leq\sum_{k=0}^{\infty}\,\sum_{n=2^{k}}^{2^{k+1}-1}\,|a_{n}|^{2}\,(1+\log_{2}n)^{2}\leq
2L<\infty.
\end{aligned}
\end{equation*}
Whence, applying the Cauchy inequality, we get
\begin{equation*}
\begin{aligned}
\sum_{k=0}^{\infty}\,\|\chi_{k}\|_{2}&=
\sum_{k=0}^{\infty}\,\|\chi_{k}\|_{2}\,(k+1)\,(k+1)^{-1}\\
&\leq\biggl(\,\sum_{k=0}^{\infty}\,\|\chi_{k}\|_{2}^{2}\,(k+1)^{2}\,\biggr)^{1/2}\;
\biggl(\,\sum_{k=0}^{\infty}\,(k+1)^{-2}\,\biggr)^{1/2}\leq\sqrt{2L}\:\sqrt{2}=
2\,\sqrt{L}.
\end{aligned}
\end{equation*}
Thus
\begin{equation}\label{eq28}
\sum_{k=0}^{\infty}\,\|\chi_{k}\|_{2}\leq2\,\sqrt{L}.
\end{equation}

Let us show that
\begin{equation}\label{eq30}
\sum_{k=0}^{\infty}\,\|\chi_{k}(x)\|_{H(x)}<\infty\quad\mbox{for}\quad\mu\mbox{-a.e}
\;\;x\in X.
\end{equation}
Recall that without loss of generality we may consider measure $\mu$ to be
$\sigma$-finite on $X$.

If $\mu(X)<\infty$, then by \eqref{eq28} and the Cauchy inequality we have:
\begin{equation}\label{eq31}
\begin{aligned}
\sum_{k=0}^{\infty}\,\int\limits_{X}\|\chi_{k}(x)\|_{H(x)}\,d\mu(x)&\leq
\sum_{k=0}^{\infty}\,\biggl(\:\int\limits_{X}\,d\mu(x)\biggr)^{1/2}
\biggl(\,\int\limits_{X}\,\|\chi_{k}(x)\|_{H(x)}^{2}\,d\mu(x)\biggr)^{1/2}\\
&\leq2\,\sqrt{\mu(X)\,L}<\infty.
\end{aligned}
\end{equation}
Therefore, according to the B.~Levi theorem, we may write
\begin{equation}\label{eq32}
\int\limits_{X}\,\biggl(\;\sum_{k=0}^{\infty}\,\|\chi_{k}(x)\|_{H(x)}\biggr)\,d\mu(x)=
\sum_{k=0}^{\infty}\:\int\limits_{X}\|\chi_{k}(x)\|_{H(x)}\,d\mu(x)<\infty;
\end{equation}
this yields \eqref{eq30}.

If $\mu(X)=\infty$, then represent $X$ as a countable union of some measurable sets
$X_{j}$, $j=1,\,2,\,3,\ldots$, with $\mu(X_{j})<\infty$. For every $j$ formula
\eqref{eq31} and its consequences, formulas \eqref{eq32} and \eqref{eq30}, remain
true if we replace $X$ by $X_{j}$. So, we get \eqref{eq30} again.

It follows from \eqref{eq30} that \eqref{eq23} is a Cauchy sequence for $\mu$-a.e.
$x\in X$, i.e., \eqref{eq23} converges. Besides,
\begin{equation*}
S^{\star}(x)\leq\sum_{k=0}^{\infty}\,\|\chi_{k}(x)\|_{H(x)}<\infty
\quad\mbox{for}\;\;\mu\mbox{-a.e.}\;\;x\in X.
\end{equation*}
Whence we have by \eqref{eq28} that:
\begin{equation}\label{eq34}
\|S^{\star}\|_{2}\leq\sum_{k=0}^{\infty}\,\|\chi_{k}\|_{2}\leq2\,\sqrt{L}.
\end{equation}

Now consider the function
\begin{equation*}
S^{\circ}(x):=\sup_{1\leq k<\infty}\,S^{\circ}_{k}(x),\quad x\in X,
\end{equation*}
where
\begin{equation*}
S^{\circ}_{k}(x):=\max_{2^{k}\leq
j<2^{k+1}}\,\bigl\|\,\sum_{n=2^{k}}^{j}\,a_{n}\,\varphi_{n}(x)\,\bigr\|_{H(x)},\quad
x\in X,\quad k=1,\,2,\,3,\ldots\,.
\end{equation*}
Applying Lemma~\ref{lem1}, with $\Psi:=(\varphi_{n})_{n=2^{k}}^{j}$ and
$b:=(a_{n})_{n=2^{k}}^{j}$, and using the condition \eqref{f3}, we may write the
following:
\begin{equation*}
\begin{aligned}
\sum_{k=1}^{\infty}\,\|S^{\circ}_{k}\|_{2}^{2}&\leq
\sum_{k=1}^{\infty}\,\max_{2^{k}\leq
j<2^{k+1}}\,\bigl(2+\log_{2}(j-2^{k}+1)\bigr)^{2}\:
\sum_{n=2^{k}}^{j}\,|a_{n}|^{2}\\&\leq
\sum_{k=1}^{\infty}\,\bigl(2+\log_{2}2^{k}\bigr)^{2}\:
\sum_{n=2^{k}}^{2^{k+1}-1}\,|a_{n}|^{2}
\leq\sum_{k=1}^{\infty}\,\sum_{n=2^{k}}^{2^{k+1}-1}\,|a_{n}|^{2}\,(2+\log_{2}n)^{2}\\
&=\sum_{n=1}^{\infty}\,|a_{n}|^{2}\,(2+\log_{2}n)^{2}\leq4L<\infty.
\end{aligned}
\end{equation*}
Therefore, by the B.~Levy theorem, we have
\begin{equation}\label{eq39}
\int\limits_{X}\,\biggl(\;\sum_{k=1}^{\infty}\,(S^{\circ}_{k}(x))^{2}\biggr)\,d\mu(x)=
\sum_{k=1}^{\infty}\,\int\limits_{X}(S^{\circ}_{k}(x))^{2}\,d\mu(x)\leq4L<\infty.
\end{equation}
Whence $\lim\limits_{k\rightarrow\infty}S^{\circ}_{k}(x)=0$ for $\mu$-a.e. $x\in X$.
This together with the convergence of \eqref{eq23} for $\mu$-a.e. $x\in X$ proved
above yields the convergence of the sequence \eqref{f3} for $\mu$-a.e. $x\in X$.

Moreover, since
\begin{equation*}
S^{*}(\Phi,a,x)\leq S^{\star}(x)+S^{\circ}(x),\quad \bigl(S^{\circ}(x)\bigr)^{2}\leq
\sum_{k=1}^{\infty}\,(S^{\circ}_{k}(x)\bigr)^{2},\quad x\in
X,
\end{equation*}
we finally deduce the required inequality \eqref{f4} from \eqref{eq34} and
\eqref{eq39},
\begin{equation*}
\|S^{*}(\Phi,a,\cdot)\|_{2}\leq\|S^{\star}\|_{2}+\|S^{\circ}\|_{2}\leq4\sqrt{L}.
\end{equation*}

Theorem \ref{th1} is proved.

\section{Proof of Theorem \ref{th2}}\label{sec5}

Without loss of generality we may assume that $a_{1}=a_{2}=\nobreak0$. Denote for an
integer $k\geq0$:
\begin{equation*}
M_{k}:=\{j\in\mathbb{N}:\,\nu_{k}+1\leq j\leq\nu_{k+1}\};
\end{equation*}
recall that $\nu_{k}:=2^{2^{k}}$. Consider an arbitrary permutation \eqref{f5} of
the orthogonal series \eqref{f1}. Define a sequence of numbers
$\bigl(\varepsilon^{(k)}_{n}\bigr)_{n=1}^{\infty}$ by the formula
\begin{equation*}
\varepsilon^{(k)}_{n}:=
  \begin{cases}
    \;1, & \text{if}\quad\sigma(n)\in M_{k}, \\
    \;0, & \text{otherwise}.
  \end{cases}
\end{equation*}
Given arbitrary $p,q\in\mathbb{N}$ with $p\leq q$, we may write
\begin{equation}\label{eq46}
\sum_{n=p}^{q}\,a_{\sigma(n)}\,\varphi_{\sigma(n)}(x)=
\sum_{k=0}^{\infty}\,\sum_{n=p}^{q}\,
\varepsilon^{(k)}_{n}\,a_{\sigma(n)}\,\varphi_{\sigma(n)}(x),\quad x\in X.
\end{equation}
The series on the right of \eqref{eq46} converges for every $x\in X$ because it
contains only a finitely many of nonzero terms.

Given any integer $k\geq0$, we set
\begin{equation}\label{eq47}
\delta_{k}(x):=\sup_{1\leq p<q<\infty}\;\bigl\|\,\sum_{n=p}^{q}\,
\varepsilon^{(k)}_{n}\,a_{\sigma(n)}\,\varphi_{\sigma(n)}(x)\,\bigr\|_{H(x)},\quad
x\in X.
\end{equation}
Note that
\begin{equation}\label{eq48}
\delta_{k}(x)\,\leq\;2\sup_{1\leq q<\infty}\;\bigl\|\,\sum_{n=1}^{q}\,
\varepsilon^{(k)}_{n}\,a_{\sigma(n)}\,\varphi_{\sigma(n)}(x)\,\bigr\|_{H(x)},\quad
x\in X;
\end{equation}
here the sum contains only the terms with $\sigma(n)\in M_{k}$. We put in,
Lemma~\ref{lem1},
\begin{equation*}
\begin{aligned}
\Psi &:=\{\varphi_{\sigma(n)}:\,
n\in\mathbb{N}\;\;\mbox{such that}\;\;\sigma(n)\in M_{k}\},\\
b &:=\{a_{\sigma(n)}:\,n\in\mathbb{N}\;\;\mbox{such that}\;\;\sigma(n)\in M_{k}\},\\
N &=N(k):=\nu_{k+1}-\nu_{k}=\nu_{k}(\nu_{k}-1).
\end{aligned}
\end{equation*}
Then
\begin{equation*}
S^{*}_{N(k)}(\Psi,b,x)=\sup_{1\leq q<\infty}\;\bigl\|\,\sum_{n=1}^{q}\,
\varepsilon^{(k)}_{n}\,a_{\sigma(n)}\,\varphi_{\sigma(n)}(x)\,\bigr\|_{H(x)},\quad
x\in X.
\end{equation*}
Therefore, by Lemma~\ref{lem1} and in view of \eqref{eq48}, we have
\begin{equation*}
\begin{aligned}
\|\delta_{k}\|_{2}&\leq(4+2\log_{2}N(k))\biggl(\;\sum_{n\,:\,\sigma(n)\in
M_{k}}\,|a_{\sigma(n)}|^{2}\biggr)^{1/2}\\
&=(4+2\log_{2}N(k))\biggl(\;\sum_{n=\nu_{k}+1}^{\nu_{k+1}}\,|a_{n}|^{2}\biggr)^{1/2}.
\end{aligned}
\end{equation*}
Hence, since
\begin{equation*}
4+2\log_{2}N(k)=4+2\log_{2}(\nu_{k}(\nu_{k}-1))\leq 8\log_{2}\nu_{k},
\end{equation*}
we arrive at the estimate
\begin{equation}\label{eq51}
\biggl(\,\int\limits_{X}\delta_{k}^{\,2}(x)\,d\mu(x)\biggr)^{1/2}\leq
8\,\biggl(\;\sum_{n=\nu_{k}+1}^{\nu_{k+1}}\,|a_{n}|^{2}\,\log^{2}_{2}n\biggr)^{1/2}.
\end{equation}

We will deduce from this that
\begin{equation}\label{eq52}
\sum_{k=0}^{\infty}\,\delta_{k}(x)<\infty \quad\mbox{для}\quad\mu\mbox{-п.в.}
\;\;x\in X.
\end{equation}
Recall that, without loss of generality, the measure $\mu$ is assumed to be
$\sigma$-finite on $X$.

If $\mu(X)<\infty$, then by the Cauchy inequality for integrals, the estimate
\eqref{eq51}, and condition \eqref{f6} we may write the following:
\begin{equation}
\begin{aligned}\label{eq53}
\sum_{k=0}^{\infty}\,\int\limits_{X}\delta_{k}(x)\,d\mu(x)&\leq
\sum_{k=0}^{\infty}\biggl(\:\int\limits_{X}d\mu(x)\biggr)^{1/2}\,
\biggl(\,\int\limits_{X}\,\delta_{k}^{\,2}(x)\,d\mu(x)\biggr)^{1/2}\\
&\leq 8\,\sqrt{\mu(X)}\;\sum_{k=0}^{\infty}\,
\biggl(\;\sum_{n=\nu_{k}+1}^{\nu_{k+1}}\,|a_{n}|^{2}\,\log^{2}_{2}n\biggr)^{1/2}<\infty.
\end{aligned}
\end{equation}
Therefore, according to the B.~Levi theorem, we have
\begin{equation}\label{eq55}
\int\limits_{X}\biggl(\;\sum_{k=0}^{\infty}\,\delta_{k}(x)\biggr)\,d\mu(x)=
\sum_{k=0}^{\infty}\:\int\limits_{X}\delta_{k}(x)\,d\mu(x)<\infty,
\end{equation}
whence we get \eqref{eq52} (recall that all $\delta_{k}\geq0$).

If $\mu(X)=\infty$, then represent $X$ as a countable union of measurable sets
$X_{j}$, $j\in\mathbb{N}$, with $\mu(X_{j})<\infty$. For every $j$ the inequality
\eqref{eq53} and its consequences, formulas \eqref{eq55} and \eqref{eq52}, remains
valid if we replace $X$ by $X_{j}$. Whence we obtain \eqref{eq52} again.

By \eqref{eq52}, for $\mu$-a.e. $x\in X$ and arbitrary $\varepsilon>0$ there exists
a number $m=m(x,\varepsilon)$ such that
\begin{equation}\label{eq56}
\sum_{k=m}^{\infty}\,\delta_{k}(x)<\varepsilon.
\end{equation}
Let $p=p(x,\varepsilon)$ be large enough so that the sum
\begin{equation*}
\sum_{n=1}^{p-1}\;a_{\sigma(n)}\,\varphi_{\sigma(n)}(x)
\end{equation*}
contains all the functions $\varphi_{n}$ whose indexes belong to $M_{k}$ with $0\leq
k< m(x,\varepsilon)$. Then by \eqref{eq47} and \eqref{eq56} we have for every $q\geq
p$  that
\begin{gather*}
\begin{aligned}
\bigl\|\,\sum_{n=p}^{q}\,a_{\sigma(n)}\,\varphi_{\sigma(n)}(x)\,\bigr\|_{H(x)}&=
\bigl\|\,\sum_{k=0}^{\infty}\,\sum_{n=p}^{q}\,
\varepsilon^{(k)}_{n}\,a_{\sigma(n)}\,\varphi_{\sigma(n)}(x)\,\bigr\|_{H(x)}\\
&=\bigl\|\,\sum_{k=m}^{\infty}\,\sum_{n=p}^{q}\,
\varepsilon^{(k)}_{n}\,a_{\sigma(n)}\,\varphi_{\sigma(n)}(x)\,\bigr\|_{H(x)}\\
&\leq\sum_{k=m}^{\infty}\;\bigl\|\,\sum_{n=p}^{q}\,
\varepsilon^{(k)}_{n}\,a_{\sigma(n)}\,\varphi_{\sigma(n)}(x)\,\bigr\|_{H(x)}\leq
&\sum_{k=m}^{\infty}\,\delta_{k}(x)<\varepsilon.
\end{aligned}
\end{gather*}

Thus, for $\mu$-a.e. $x\in X$ and for an arbitrary $\varepsilon>0$ there exists a
number $p=p(x,\varepsilon)$ such that
\begin{equation*}
\bigl\|\,\sum_{n=p}^{q}\,a_{\sigma(n)}\,\varphi_{\sigma(n)}(x)\,\bigr\|_{H(x)}
<\varepsilon
\end{equation*}
for every integer $q\geq p$. So, the series \eqref{f5} converges for $\mu$-a.e.
$x\in X$.

Theorem~\ref{th2} is proved.

\section{Proof of Theorem~\ref{th3}}\label{sec6}

We deduce it from Theorem~\ref{th2} by showing that the conditions \eqref{f7} and
\eqref{f8} together imply~\eqref{f6}.

For every integer $k\geq0$, put
\begin{equation*}
A_{k}:=\sum_{n=\nu_{k}+1}^{\nu_{k+1}}\, |a_{n}|^{2}\,\log_{2}^{2}n;
\end{equation*}
here $\nu_{k}:=2^{2^{k}}$ as above. Applying the Cauchy inequality, we may write
\begin{equation*}
\sum_{k=0}^{\infty}\,A_{k}^{1/2}=
\sum_{k=0}^{\infty}\,A_{k}^{1/2}\,\omega_{\nu_{k}}^{1/2}\,\omega_{\nu_{k}}^{-1/2}\leq
\biggl(\,\sum_{k=0}^{\infty}\,A_{k}\,\omega_{\nu_{k}}\biggr)^{1/2}\;
\biggl(\,\sum_{k=0}^{\infty}\,\omega_{\nu_{k}}^{-1}\biggr)^{1/2}.
\end{equation*}
It is known that
\begin{equation*}
\sum_{n=2}^{\infty}\;\frac{1}{n\,(\log_{2}n)\,\omega_{n}}<\infty\;\Leftrightarrow\;
\sum_{n=1}^{\infty}\;\frac{1}{n\,\omega_{2^{n}}}<\infty\;\Leftrightarrow\;
c:=\sum_{n=0}^{\infty}\;\frac{1}{\omega_{\nu_{n}}}<\infty.
\end{equation*}
Therefore, using \eqref{f7} and since $(\omega_{n})_{n=1}^{\infty}$ is increasing, we
have the following:
\begin{equation*}
\begin{aligned}
\biggl(\,\sum_{k=0}^{\infty}\,A_{k}^{1/2}\biggr)^{2}&\leq
c\,\sum_{k=0}^{\infty}\,A_{k}\,\omega_{\nu_{k}}=
c\,\sum_{k=0}^{\infty}\,\omega_{\nu_{k}}\,\sum_{n=\nu_{k}+1}^{\nu_{k+1}}\,
|a_{n}|^{2}\,\log_{2}^{2}\\
&\leq c\,\sum_{k=0}^{\infty}\,\sum_{n=\nu_{k}+1}^{\nu_{k+1}}\,
|a_{n}|^{2}\,(\log_{2}^{2}n)\,\omega_{n}=
c\,\sum_{n=3}^{\infty}\,|a_{n}|^{2}\,(\log_{2}^{2}n)\,\omega_{n}<\infty.\label{eq65}
\end{aligned}
\end{equation*}
Thus, the condition \eqref{f6} is satisfied,
\begin{equation*}
\sum_{k=0}^{\infty}\,\biggl(\;\sum_{n=\nu_{k}+1}^{\nu_{k+1}}\,
|a_{n}|^{2}\,\log_{2}^{2}n\biggr)^{1/2}=\sum_{k=0}^{\infty}\,A_{k}^{1/2} <\infty.
\end{equation*}
Therefore, by Theorem~\ref{th2}, the sequence \eqref{f1} converges unconditionally
$\mu$-a.e. on $X$.

Theorem~\ref{th3} is proved.

\section{Final remark}\label{sec7}

A simple inspection of the proofs of Lemma~\ref{th1} and Theorems
\ref{th1}--\ref{th3} reveals that they remain true if the system
$(\varphi_{n})_{n=1}^{\infty}$ forms a Riesz basis in the closure of its linear span
in $\mathbf{L}_{2}$. In this case, the factor $C\log_{2}(N+1)$ should be used, instead of
$2+\nobreak\log_{2}N$, in the right-hand side of \eqref{eq10}, the constant $C>0$ as
well as $K$ in Theorem~\ref{th1} depending on a choice of
$(\varphi_{n})_{n=1}^{\infty}$.

\end{document}